\documentclass[subeqn, letterpaper, reqno, 11pt]{amsart}
\usepackage[T1]{fontenc}
\usepackage{amssymb, amsmath, amsthm, graphicx, mathtools, dsfont, doi}
\usepackage{setspace} % line spacing
\usepackage[numbers, sort&compress]{natbib}
\usepackage[all]{xy}
\usepackage{hyperref}
\usepackage[fulloldstylenums]{kpfonts}
\DeclareMathAlphabet\mathfrak{U}{euf}{m}{n}
\SetMathAlphabet\mathfrak{bold}{U}{euf}{b}{n}

\newcommand{\R}{\mathbb{R}}
\def\C{\mathbb{C}}
\newcommand{\T}{\mathbb{T}}
\newcommand{\sS}{\mathbb{S}}
\newcommand{\Ci}[1]{\mathcal{C}^{#1}}
\newcommand{\Par}{\mathcal{P}}
\newcommand{\Hol}{\mathcal{H}}
\newcommand{\OC}{\mathcal{OC}}
\newcommand{\X}{\frak{X}}
\newcommand{\inv}{^{\text{-}1}}
\newcommand{\n}{^{\natural}}
\newcommand{\lie}{\frak{g}}

\newcommand{\zero}{\textup{\fontfamily{jkp}\selectfont 0}}

\swapnumbers %Numbers before def, props, thms, etc.
\numberwithin{equation}{section}

\newtheorem{theorem}{Theorem}[section]
\newtheorem{lem}[theorem]{Lemma}
\newtheorem*{thma}{Theorem A}
\newtheorem*{thmb}{Theorem B}
\newtheorem*{thmc}{Theorem C}

\newtheorem{proposition}[theorem]{Proposition}
\newtheorem{cor}[theorem]{Corollary}
\newtheorem{co}[theorem]{Conjecture}

\newtheorem{pro}{Problem}
\newtheorem{que}[theorem]{Question}

\theoremstyle{definition}
\newtheorem{df}[theorem]{Definition}
\newtheorem{exa}[theorem]{Example}

\theoremstyle{remark}
\newtheorem{rem}[theorem]{Remark}
\newtheorem*{cla}{Claim}

\newcommand{\bd}{\begin{df}}
\newcommand{\ed}{\end{df}}

\newcommand{\bq}{\begin{que}}
\newcommand{\eq}{\end{que}}

\newcommand{\bl}{\begin{lem}}
\newcommand{\el}{\end{lem}}

\newcommand{\br}{\begin{rem}}
\newcommand{\er}{\end{rem}}

\newcommand{\bt}{\begin{theorem}}
\newcommand{\et}{\end{theorem}}

\newcommand{\bc}{\begin{cor}}
\newcommand{\ec}{\end{cor}}

\newcommand{\bco}{\begin{co}}
\newcommand{\eco}{\end{co}}

\newcommand{\bcl}{\begin{cla}}
\newcommand{\ecl}{\end{cla}}

\newcommand{\bp}{\begin{proposition}}
\newcommand{\ep}{\end{proposition}}

\newcommand{\brm}{\begin{rem}}
\newcommand{\erm}{\end{rem}}

\newcommand{\be}{\begin{equation}}
\newcommand{\ee}{\end{equation}}
\newcommand{\bx}{\begin{exa}}
\newcommand{\ex}{\end{exa}}
\newcommand{\bpr}{\begin{pro}}
\newcommand{\epr}{\end{pro}}

\begin{document}

\title{Waning in principal bundles}
\author{Pedro Sol\'{o}rzano}\thanks{The author is supported by IMPA and CAPES}
\address{Departamento de Matem\'{a}tica\\
Universidade Federal de Santa Catarina\\
Campus Universit\'{a}rio Trindade\\
CEP 88.040-900 Florian\'{o}polis-SC, Brasil}
\email{solorzano@member.ams.edu}
\date{March 31, 2015}

\maketitle
\begin{abstract} 
Let $P\to M$ be a principal bundle. Consider a sequence of metrics on $P$ obtained by re-scaling the fibers to points. The Gromov-Hausdorff limit of the tangent bundles over these principal bundles with their Sasaki metric is seen herein to be a locally trivial fiber bundle containing the tangent space to the base as a subbundle in a natural way.  Berger 3-spheres provide an example where the limit fibers are still of dimension 3. The fibers are shown to be entirely determined by the riemannian holonomy of the chosen bi-invariant metric.
\end{abstract}

\section*{Introduction}
Understanding the behavior of geometric structures as they undergo a limiting process has remained an important topic since its introduction in the late 1970's by M. Gromov. It has been instrumental in the proofs of big conjectures as recent as the work of \citet{MR3261011}.  In particular, there they have to deal with the convergence of powers of line bundles.  In this communication, the purpose is to study a particular example of vector bundles---namely tangent bundles---over principal bundles. 

Processes that collapse fibers to a point abound. R. Hamilton's much celebrated Ricci Flow has also provided examples and tools in this direction. In the case of the K\"ahler-Ricci flow \citet*{MR3010688} show that a projective bundle over a projective base collapses the fibers as it undergoes this flow in finite time if the initial metric is in a suitable K\"ahler class. 

Restrict the attention to vector or to sphere bundles. In the compact case, for unit tangent (sub)bundles, \citet*{MR2521897} studied the behavior of a sequence of horizontal unit sphere bundles, with the requirement of integrability and under a general warping. In particular they show that these horizontal unit sphere bundles converge to the tangent sphere bundle of the base if and only if the total space converges to the base.  They analyze metrics of type $(p,q)$ introduced by \citet*{MR2330461}.  This report deals with a similar result, though not a direct generalization since it deals only with the Sasaki metric (which is the case $p=q=0$), removes both the integrability and the horizontality restriction, and focuses on the setting of principal bundles. 

In \cite{MR3261732}, the author verified the pointed Gromov-Hausdorff convergence of the Sasaki metrics on the tangent bundles (both the total spaces and their projection maps, up to passing to a subsequence) over a sequence of converging riemannian spaces in the Gromov-Hausdorff sense. Furthermore, the structure of the limit is such that the fiber of the limit map at a point $p$ is homeomorphic to a quotient of the form
\be
\R^k/G
\ee
where G is called {\em wane group} at $p$ of the sequence. The wane group is well-defined up to conjugation by an orthogonal matrix and depends on the point and the sequence. The purpose of this communication is to study the case of connection metrics on principal $G$ bundles, collapsing to the base space by re-scaling the fibers. In particular, in this case the wane groups are all equal to the closure of the holonomy group of the bi-invariant metric on $G$.
\begin{thma}[Homogeneous Waning Theorem]\label{HWT} Let $P$ be a principal $G$-bundle with connection $\omega$ over a riemannian manifold $(M,g_M)$ with projection map $\pi:P\rightarrow M$. Let $\{\varepsilon_n\}$ be a sequence of positive numbers converging to zero. Let $g_n$ be the connection metric on $P$ given by 
\be
g_n=\pi^*g_M+\varepsilon_n^2g_G
\ee
for a fixed bi-invariant metric $g_G$ on $G$. Let $g_n^S$ be Sasaki metrics on $TP$ induced from $g_n$. Then, for any point $x\in M$ and for any point $p\in\pi\inv(x)$, the sequence of pointed spaces $(TP,g_n^S,\zero_p)$ with their projections $\pi_P$ converges in the pointed Gromov-Hausdorff sense to a space $Y$ with projection $\pi_Y:Y\rightarrow X$ with wane groups of $G_x$ given by
$$
\overline{Hol(g_G)}.
$$
\end{thma}
The proof of this statement is the content of Section \ref{proof}. In particular, given a sequence of principal $\T^n$-bundles over a base manifold $(M,g)$ with connection metrics and collapsing fibers, the corresponding wane groups are all equal and trivial, that is  $G=\{I_k\}$ regardless of the base point. This is the case for example of the Berger spheres converging to a 2-sphere. 

Topologically, this is the first non trivial example of the survival of the tangent spaces even when the spaces collapses in dimension. Furthermore, these maps $\pi_Y$ are locally trivial.
\begin{thmb} Let $P$ be a principal $G$-bundle with connection $\omega$ over a riemannian manifold $(M,g_M)$ with projection map $\pi:P\rightarrow M$. Let ${\epsilon_n}$ be a sequence of positive numbers converging to zero. Let $g_n$ be the connection metric on $P$ given by 
\be
g_n=\pi^*g_M+\varepsilon_n^2g_G
\ee
for a fixed bi-invariant metric $g_G$ on $G$. Let $g_n^S$ be Sasaki metrics on $TP$ induced from $g_n$. Then, the pointed Gromov-Hausdorff limit $Y$ of the sequence of pointed spaces $TP$ is a locally trivial bundle that contains the tangent bundle $TM$ in a natural way.
\end{thmb}
The proof of this is the content of Section \ref{prooflt}. It should be noted that no claim is made with respect to the structure group. This in particular describes the structure of the limit of the Sasaki metrics on the Berger 3-spheres.

In fact, the proof of Theorem B provided requires the analysis of parallelism for the limits. In \cite{MR3261732}, the author proved that there is a weak notion of parallel translation along curves in the limit spaces. This is accomplished by saying that a curve is horizontal in $Y$ if it projects via $\pi_Y$ to a curve of the same length in $M$. Existence of horizontal lifts given initial conditions is proved in \cite{MR3261732}.
\begin{thmc} Let $P$ be a principal $G$-bundle with connection $\omega$ over a riemannian manifold $(M,g_M)$ with projection map $\pi:P\rightarrow M$. Let ${\epsilon_n}$ be a sequence of positive numbers converging to zero. Let $g_n$ be the connection metric on $P$ given by 
\be
g_n=\pi^*g_M+\varepsilon_n^2g_G
\ee
for a fixed bi-invariant metric $g_G$ on $G$. Let $g_n^S$ be Sasaki metrics on $TP$ induced from $g_n$.  Then, the pointed Gromov-Hausdorff limit $Y$ of the sequence of pointed spaces $TP$ parallel translation in $Y$ along curves in $M$ is unique.
\end{thmc}

This is the content of section \ref{unipar}. The author had also proved that uniqueness implies that all wane groups must coincide up to conjugation.  Theorem C provides evidence of the converse.

Recent work of \citet*{1502.05307} suggests that even in the more general case of principal actions, the conclusion might still hold replacing $Hol(G)$ with $Hol(G/H)$ and a constant re-scaling with a Cheeger deformation. 

\

%%%%%%%%%%%%%%%%%%%%%%%%%%%%%%%%%%%%%%%%%%%%%%%
%%%%%%%%%%%%%%%%%%%%%%%%%%%%%%%%%%%%%%%%%%%%%%%
% A C K N O W L E D G E M E N T S
%%%%%%%%%%%%%%%%%%%%%%%%%%%%%%%%%%%%%%%%%%%%%%%
%%%%%%%%%%%%%%%%%%%%%%%%%%%%%%%%%%%%%%%%%%%%%%%

{\bf\noindent Acknowledgements.} The author is thankful of the support given by IMPA and CAPES through their postdoctorate of excellence fellowships. The author wishes to thank P. Gianniotis, R. Perales and C. Pro for meaningful discussions that lead to this report. In particular R. Perales for reading early drafts and for her invaluable suggestions. 

%%%%%%%%%%%%%%%%%%%%%%%%%%%%%%%%%%%%%%%%%%%%%%%
%%%%%%%%%%%%%%%%%%%%%%%%%%%%%%%%%%%%%%%%%%%%%%%

%%%%%%%%%%%%%%%%%%%%%%%%%%%%%%%%%%%%%%%%%%%%%%%
%%%%%%%%%%%%%%%%%%%%%%%%%%%%%%%%%%%%%%%%%%%%%%%

\setcounter{tocdepth}{2}
\tableofcontents
\pagebreak

%%%%%%%%%%%%%%%%%%%%%%%%%%%%%%%%%%%%%%%%%%%%%%%
%%%%%%%%%%%%%%%%%%%%%%%%%%%%%%%%%%%%%%%%%%%%%%%

\section{Background}
In this section an assortment of results that will be needed in the sequel are stated. The first subsection deals with writing parallel translation in local form as a matrix solution for an ODE.  The second is a description of the Saki-type metrics. The third reviews the definitions of the Gromov-Hausdorff theory of convergence together with a suitable version of the Arzel\`a-Ascoli Theorem.  The fourth subsection reviews the notion of waning, central to this communication as well as some results about parallelism. The last section introduces notation from principal bundles and recalls the basic formulas for covariant differentiation for connection metrics.  
\subsection{Lemmata about parallel translation}
The following two results are stated here for completeness and for future reference. These are elementary and references are given. 
\bl\label{parallim}  Given a sequence of curves $\alpha_n$ converging in the $\Ci{1}$ sense to a curve $\alpha$, their corresponding parallel ${P}^{\alpha_n}$ translations converge to ${P}^{\alpha}$.
\el
\begin{proof}
It directly follows from an exercise in the book by \citet{MR2229062} once one observes that for large enough $n$ all curves $\alpha_n$ are homotopic to $\alpha$.
\end{proof}

Parallel translation along a curve $\alpha$ can be described as an ordinary differential equation once a local orthonormal frame is given. This can be written as in matrix form as
\be
(\dot P^i)=(Q_j^i)\cdot(P^i),
\ee
where $(Q_j^i)\in\mathfrak{so}(n)$, given the representation $\nabla_{\dot\alpha}E_i=Q_i^jE_j$.
In particular, one has the following general lemma about this type of linear ODEs. 

\bl \label{paramdepd}  Consider sequences $\{Y_n(t)\}$ and $\{Q_n(t)\}$ such that:
\begin{enumerate} 
\item they satisfy the equation $Y'_n=Q_n\cdot Y_n$  for $t\in [\zero,1]$ with the same initial condition $Y_n(\zero)=Y_\zero$ for all n;  
\item there exists a uniform $K>\zero$ such that $\|Y_n(t)\|\leq K$.
\item the sequence $\{Y_n(1)\}$ is convergent and it converges to $Y_1$; and that
\item  there is a linear subspace $\lie$ of the space of matrices such that the sequence $\{t\mapsto dist(Q_n(t),\lie)\}$ is converging uniformly to zero on $[\zero,1]$.
\end{enumerate}
Then  $Y_1$ is also a limit of a sequence  $\{Z_n(1)\}$, where 
\begin{enumerate}
\item The sequence $\{Z_n\}$ is a sequence of solutions to  $Z'_n=L_n\cdot Z_n$  with $Z_n(\zero)=Y_\zero$; and 
\item for any $t\in[\zero, 1]$, $L_n(t)\in \lie$.
\end{enumerate}
\el
\proof The case when the sequence $\{Q_n\}$ is itself bounded follows from the general theory of smooth dependance on initial conditions and parameters. For example, by regarding the problem as a matrix problem,
\be
\left\{
\begin{array}{c} 
\dot Y=Q\cdot Y\\
Y(\zero)=I
\end{array}\right.
\ee
for small enough operator norm, the solutions can always be expressed as 
\be
Y=\exp\Omega,
\ee
where $\Omega$ is an integral expression akin to the Baker-Campbell-Hausdorff formula for $Q$, $\int_\zero^tQ(\tau)\,d\tau$ and commutators thereof, as suggested by \citet{MR0067873}. 

Even if the sequence is unbounded, let $L_n$ be the projection of $Q_n$ unto $\lie$ and let $\{Z_n\}$ satisfy 
$$Z'_n=Q_n\cdot Z_n,$$
with same initial condition $Z_n(\zero)=Y_\zero$. Let $\varepsilon_n=(Q_n-L_n)Y_n$. Observe that
\begin{enumerate}
\item The norm $\|\varepsilon_n\|$ goes uniformly to zero as $n$ goes to infinity. 
\item The curves $X_n:=Y_n-Z_n$ satisfy the differential equation
$$
X'_n=L_n\cdot X_n+\varepsilon_n,
$$
with initial condition $X_n(\zero)=\zero$. 
\end{enumerate}
Whence $\{X_n\}$ converges up to passing to a subsequence to the constant solution. 
\endproof

\br In particular one has that should $\|(Q_i^j)\|$ be small, then the frame is close to being parallel. 
\er

\subsection{Sasaki Metrics}

The Sasaki metric is described in terms of a slightly more general metric, the Sasaki-type metric on vector bundles, as follows. Given a euclidean metric $h$ with compatible connection $\nabla$ on a vector bundle over a manifold with riemannian metric $g$, the corresponding Sasaki-type metric is the unique metric such that for any curve $\gamma$ in the total space, with projected curve $\delta$,
\be\label{sasdefeq}
\|\dot\gamma\|^2=\|\dot\delta\|_g^2+\|\nabla_{\dot\delta}\gamma\|_h^2.
\ee
On the tangent bundle one immediately has a Sasaki metric associated to the metric $h=g$ and the Levi-Civita connection.

\subsection{Gromov-Hausdorff convergence}
Listed are a few results needed in the sequel. For a detailed treatment of this concepts see \cite{MR1835418,MR2307192, MR2243772}. 

\bd[\citet{MR2307192}] Given two complete metric spaces $(X,d_X)$ and $(Y,d_Y)$, their {\em Gromov-Hausdorff} distance is defined as the following infimum.
\be\label{GHdist}
d_{GH}(X,Y)=\inf\left\{\varepsilon>\zero\left|\begin{array}{rl}(1)&\exists d:(X\sqcup Y)\times(X\sqcup Y)\rightarrow\R,\text{ metric}\\ (2)& d|_{X\times X}=d_X,d|_{Y\times Y}=d_Y \\ (3)&\forall x\in X(\exists y\in Y, d(x,y)<\varepsilon) \\ (4)&\forall y\in Y(\exists x\in X, d(x,y)<\varepsilon)\end{array}\right. \right\}
\ee
That is the infimum of possible $\varepsilon>\zero$ for which there exists a metric on the disjoint union $X\sqcup Y$ that extends the metrics on $X$ and $Y$, in such a way that any point of $X$ is $\varepsilon$-close to some point of $Y$ and vice versa.  
\ed
\bd[\citet{MR2307192}] Let $X$ and $Y$ be metric spaces. For $\varepsilon>\zero$, an {\em $\varepsilon$-isometry} from $X$ to $Y$ is a (possibly non-continuous) function $f:X\rightarrow Y$ such that:
\begin{enumerate}
\item for all $x_1,x_2\in X$,
\be
\left|d_X(x_1,x_2)-d_Y(f(x_1),f(x_2))\right|<\varepsilon;\text{ and}
\ee 
\item for all $y\in Y$ there exists $x\in X$ such that
\be 
d_Y(f(x),y)<\varepsilon.
\ee
\end{enumerate}
\ed

\bp[\citet{MR2307192}]\label{epsiso} Let $X$ and $Y$ be metric spaces and $\varepsilon>\zero$. Then,
\begin{enumerate}
\item if $d_{GH}(X,Y)<\varepsilon$ then there exists a $2\varepsilon$-isometry between them. 
\item if there exists an $\varepsilon$-isometry form $X$ to $Y$, then $d_{GH}(X,Y)<2\varepsilon$.
\end{enumerate}
\ep
\br\label{GHisH} It was proved by \citet{MR2307192} that if a sequence $\{X_i\}$ of compact metric spaces converges in the Gromov-Hausdorff sense to a compact metric space $X$, then there exists a metric on $X\sqcup\bigsqcup_iX_i$ for which the sequence $\{X_i\}$ converges in the Hausdorff sense. Because of this, it makes sense to say that a sequence of points $x_i\in X_i$ converge to a point $x\in X$. 

\er

\br\label{ptsqnHaus} In this setting, a sequence of subspaces $X_i\subseteq Z$ converges to a subspace $X\subseteq Z$ if and only if: 
\begin{enumerate}
\item for any convergent sequence $x_i\rightarrow x$, such that $x_i\in X_i$ for all $i$, it follows that $x\in X$; and
\item for any $x\in X$ there exists a convergent sequence $x_i\rightarrow x$, with $x_i\in X_i$.
\end{enumerate}

\er

\bd[\citet{MR2307192}]\label{ptsqns} Let $\{X_i\}$, $\{Y_i\}$ be convergent sequences of compact pointed metric spaces and let $X$ and $Y$ be their corresponding limits. One says that a sequence of continuous functions $\{f_i\}:\{X_i\}\rightarrow\{Y_i\}$ {\em converges} to a function $f:X\rightarrow Y$ if there exists a metric on $X\sqcup\bigsqcup_iX_i$ for which the subspaces $X_i$ converge in the Hausdorff sense to $X$ and such that for any sequence $\{x_i\in X_i\}$ that converges to a point $x\in X$, the following holds.
\be 
f(x)=\lim_{i\rightarrow\infty}f_i(x_i)
\ee
\ed
\br The limit function $f$ is unique if it exists; i.e. it is independent of the choice of metric on $X\sqcup\bigsqcup_iX_i$.
\er

The following is Gromov's way to produce a notion of convergence for the non-compact case. For technical reasons, the assumption that the spaces be proper (i.e. that the distance function from a point is proper, thus yielding that closed metric balls are compact) is required \cite{MR2307192}. 

\bd \label{ptGH}A sequence $\{(X_i,x_i)\}$ of pointed proper metric spaces is said to converge to $(X,x)$ in the pointed Gromov-Hausdorff sense if the following holds: For all $R>\zero$ and for all $\varepsilon>\zero$ there exists $N$ such that for all $i>N$ there exists an $\varepsilon$-isometry 
$$f_i:B_R(x_i)\rightarrow B_R(x),$$
 with $f_i(x_i)=x$, where the balls are endowed with restricted (not induced) metrics. 
\ed
In the previous definition, it is enough to verify the convergence on a sequence of balls around $\{x_i\}$ such that their radii $\{R_j\}$ go to infinity. Furthermore, the limit is necessarily also proper as noted by \citet{MR623534}. Given a pointed space $(X,x)$, \citet{MR2307192} studies the relation between precompactness and the function that assigns to each choice of $R>\zero$ and $\varepsilon>\zero$ the maximum number $N=N(\varepsilon, R,X)$ of disjoint balls of radius $\varepsilon$ that fit within the ball of radius $R$ centered at an $x \in X$. Furthermore he proves the following result.

\bt[Gromov's Compactness Theorem \cite{MR2307192}, Prop. 5.2 ]\label{Gromovs} A family of pointed path metric spaces $(X_i,x_i)$ is pre-compact with respect to the pointed Gromov-Hausdorff convergence if and only if each function $N(\varepsilon, R, \cdot)$ is bounded on $\{X_i\}$. In this case, the family is relatively compact, i.e., each sequence in the $X_i$ admits  a subsequence that converges in the pointed Gromov-Hausdorff sense to a complete, proper path metric space. 
\et

\br \label{GromovN}Providing a bound for $N$ is equivalent to providing a bound to the minimum number of balls of radius $2\varepsilon$ required to cover the ball of radius $R$ (see \cite{MR2243772}). This will be used instead in the sequel and will also be denoted by $N$.
\er

\br In the non-compact setting, in order to consider the convergence of sequences of points $\{p_i\in X_i\}$, as in Remark \ref{GHisH}, the only technicality is the following: In order for a sequence $\{p_i\}$ to be convergent, it has to be bounded. Therefore, there must exist a large enough $R>\zero$ such that for all $i$, $p_i\in B_R(x_i)$, where the $x_i\in X_i$ are the distinguished points. Because of this, a sequence  $\{p_i\in X_i\}$ is convergent if there exists $R>\zero$ for which the sequence  $\{p_i\in \overline{B_R(x_i)}\subseteq X_i\}$ is convergent as in Remark \ref{GHisH}.
\er
\bt[Arzel\`a-Ascoli Theorem, cf. \cite{MR623534, MR1094462,MR3261732}]\label{ArAs}
Consider two convergent sequences of complete proper pointed metric spaces, $\{(X_i,x_i)\}$, $\{(Y_i,y_i)\}$. Let $(X,x)$ and $(Y,y)$ be their corresponding limits. Suppose further that there is an equicontinuous sequence of continuous maps $\{f_i\}$,
\be 
f_i:X_i\rightarrow Y_i,
\ee
such that $f_i(x_i)=y_i$, for all $i$. Then there exist a continuous function $f:X\rightarrow Y$, with $f(x)=y$, and a subsequence of $\{f_i\}$ that converges to $f$. 
\et

\subsection{Wane groups and the convergence of Sasaki-type metrics}

This subsection is a quick review of the results proved by the author in \cite{MR3261732}. They provide both a framework and the motivation for this report. 

\bt \label{precptBWC} Given a precompact collection of (pointed) riemannian manifolds $\mathcal{M}$ and a positive integer $k$, the collection $\textup{BWC}_k(\mathcal{M})$ of vector bundles with metric connections of rank $\leq k$ endowed with metrics of Sasaki-type is also precompact. The distinguished point for each such bundle is the zero section over the distinguished point of their base.
\et
In particular, for a convergent sequence, a description of the limit in the vertical direction is given.

\bt\label{finholtype} Let $\pi_i:E_i\rightarrow X_i$ be a convergent sequence of vector bundles with bundle metric and compatible connections $\{(E_i,h_i, \nabla_i)\}$, with limit $\pi:E\rightarrow X$. Then there exists a positive integer $k$ such that for any point $p\in X$ there exists a compact Lie group $G\leq O(k)$, that depends on the point, such that the fiber $\pi\inv(p)$ is homeomorphic to $\R^k/G$, i.e. the orbit space under the standard action of $G$ on $\R^k$. This group is called the {\em wane group} of the sequence $\{\pi_i\}$.
\et

A description of these wane groups was also given. In particular, one has to look at the structure of the fibers as {\em holonomic spaces}.

\bd[\citet{1004.1609}] A {\em holonomic space} is a triplet $(V,H,L)$: $V$ is a normed vector space, $H$ is a group of norm preserving linear maps, and $L$ is a group-norm on $H$; such that together satisfy that for any $u \in V$ there exists a positive number $R>\zero$ such that if $v,w\in V$ and for any $a\in H$, if
$$\|u-v\|,\|u-w\|<R$$
then
\be\label{propP}
\|v-w\|^2-\|av-w\|^2\leq L^2(a).
\ee
\ed

\bt[\citet{1004.1609}] \label{hlmetthm} Let $(V,H,L)$ be a holonomic space.  
\be\label{hlmet} d_L(u,v)=\inf_{a\in H}\left\{\sqrt{L^2(a)+\|u-av\|^2}\right\},
\ee
is a metric on $V$. Furthermore, the identity map is a distance non-increasing map between $(V,\|\cdot\|)$ and $(V,d_L)$.
\et

\bd[\cite{1004.1609}]  The metric given by \eqref{hlmet} is called {\em holonomic space metric}.
\ed

In \cite{1004.1609}, the author proves that the restricted metric on the fibers of vector bundle are indeed holonomic spaces; namely, the following results.

\bt[\citet{1004.1609}]\label{lengthproof} Let $Hol_p$ be the holonomy group over a point $p\in M$ of a bundle with metric and connection and suppose that $M$ is riemannian. Then the function $L_p:Hol_p\rightarrow \R$, 

\be\label{lengthnorm}
L_p(A)=\inf\{\ell(\alpha)|\alpha\in\Omega_p,P^{\alpha}_1=A\},
\ee
is a group-norm for $Hol_p$. Here $\Omega_p$ denotes the space of piece-wise smooth loops based at $p$. 
\et
\bd[\citet{1004.1609}]\label{lengthdef}
The function $L_p$, defined by \eqref{lengthnorm} is called {\em length norm} of the holonomy group induced by the riemannian metric at $p$.
\ed

\bt[\citet{1004.1609}]\label{holfib} Let $E_p$  be the fiber of a vector bundle with metric and connection $E$ over a riemannian manifold $M$ at a point $p$. Let $Hol_p$ denote the associated holonomy group at $p$ and let $L_p$ be the group-norm given by \eqref{lengthnorm}. Then $(E_p,Hol_p, L_p)$ is a holonomic space. Moreover, if $E$ is endowed with the corresponding Sasaki-type metric, the associated holonomic distance coincides with the restricted metric on $E_p$ from $E$.
\et

And then, the following result gives a more precise description of the wane groups of Theorem \ref{finholtype}.

\bt[\citet{MR3261732}]\label{limsup} Let $V$ be a finite dimensional normed vector space and $\{H_i\}$ be a sequences of subgroups of the group of norm preseving linear maps, denoted here by $O(V)$. Consider a sequence of group-norms $\{L_i:H_i\rightarrow\R\}$ such that the semi-metrics $d_i=d_{L_i}$, given by
\be\label{defholsem}
d_{L_i}(u,v)=\inf_a\sqrt{L_i^2(a)+\|au-v\|^2},
\ee
for any $u, v \in V$, converge uniformly on compact sets to a semi-metric $d_{\infty}$ on $V$. Then, there exists a closed subgroup $G\subseteq O(V)$ given by
$$
G=\OC_{O(V)}(\langle G_\zero\rangle),
$$
the orbit closure\footnote{The orbit closure in $O(V)$ of a subgroup $A\leq O(V)$ is the largest subgroup $B\in O(V)$ having the same orbits as $A$.} in $O(V)$ of the group generated by the set
\be
G_\zero=\{g\in O(V)| g=\lim_{i_n\rightarrow\infty}a_{i_n}, \lim_{i_n\rightarrow\infty}L_{i_n}(a_{i_n})=\zero\},
\ee
such that for any $u,v\in V$,  
\be 
d_{\infty}(u,v)=\zero
\ee
if and only if there exists $g\in G$ such that $v=gu$.
\et

Another important fact about these limits is that they have a weak notion of parallelism. The following can also be found in \citet{MR3261732}.

\bd Given a surjective submetry $\pi:Y\rightarrow X$ a curve $\gamma:[\zero,1]\rightarrow Y$ is {\em horizontal} if and only if
\be
\ell(\gamma)=\ell(\pi\gamma).
\ee
The set of all such curves will be denoted by $\Hol(\pi)$.
\ed

\bd Given a curve $\alpha:[\zero,1]\rightarrow X$ and a point $u\in\pi\inv\alpha(\zero)$ a {\em parallel translation of $u$ along $\alpha$} is a horizontal $\gamma$ such that $\gamma(\zero)=u$ and $\pi\gamma=\alpha$.
\ed

\bd Given a curve $\alpha:[\zero,1]\rightarrow X$ and a point $u\in\pi\inv\alpha(\zero)$ {\em the parallel translation of $u$ along $\alpha$} is given as relation $\Par\subseteq \pi\inv(\alpha(\zero))\times\pi\inv(\alpha(1))$. This can be regarded as a set-valued function 
$$ \Par^{\alpha}:\pi\inv(\alpha(\zero))\dashrightarrow\pi\inv(\alpha(1)),
$$
given by
\be \Par^{\alpha}(u)=\{\gamma(1)| \gamma\in\Hol(\pi), \pi\gamma=\alpha\}\subseteq\pi\inv(\alpha(1)).
\ee
\ed
And in particular,
\bt  For any sequence of riemannian manifolds $\{(X_i,p_i)\}$ converging to $(X,p)$ consider a convergent family of bundles with metric connection $(E_i,h_i,\nabla_i)$ over it converging to $(E,\varsigma(p))$  with $\pi:E\rightarrow X$. Let $\alpha:I\rightarrow X$ be any rectifiable curve. Then for any $u\in \pi\inv(\alpha(\zero))$ there exists  $v\in \pi\inv(\alpha(1))$ such that there exists a horizontal path $\gamma$ from $u$ to $v$. In other words, 
\be
\Par^{\alpha}(u)\neq\varnothing.
\ee
\et
Furthermore,  the following result which will be used in the sequel. 
\bt\label{limhh} Horizontal curves are the uniform limits of horizontal curves.
\et

\br\label{absnorm} Since one also has that the norms $\|TM\|\rightarrow\R$ converges to a map $\mu$, horizontal curves have constant norm and constant re-scalings of horizontal curves are horizontal. (See \citet{MR3261732}).
\er

\subsection{Principal Bundles and Connection Metrics}
Most of the results in this section are fairly classical (see the book by \citet{MR1393940}). The purpose of this subsection is to introduce the notation and to state Proposition \ref{covdevP}, which is instrumental in the sequel. 

A {\em principal $G$-bundle} $\xi=(P,M, \pi:P\rightarrow M, G, \mu:P\times G\rightarrow P)$ consists of a total space $P$, a base space $M$, a Lie group $G$ together with a free right action $\mu$ such that the projection $\pi$ is equivalent to the canonical projection from $P$ to $P/G$. 

Let $\lie$ be the Lie algebra of $G$. Then a {\em connection} is a $\lie$-valued 1-form $\omega$ such that
\be
(R_a)^*\omega= \textup{ad}(a\inv)\omega
\ee
and
\be
\omega(A^*)=A
\ee
where $R_a(p)=\mu(p,a)$ and $A^*$ is the image of $A\in\lie$ under the canonical Lie algebra homomorphism between $\lie$ and $\X(P)$. 

With this data, a {\em connection metric} is a riemannian metric on $P$ obtained from a riemannian metric $\langle\ , \ \rangle_M$ on $M$ and an invariant metric $\langle\ , \ \rangle_G$ on $G$, in the following way:

\be
\langle X,Y  \rangle_P=\langle\pi_*X , \pi_*Y \rangle_M+\langle\omega(X) ,\omega(Y) \rangle_G
\ee
This is clearly a riemannian metric on P since the horizontal distribution 
\be
Q_p=\{X_p\in T_pP|\ \omega(X)=\zero \}
\ee
is complementary to the vertical  distribution $G_p=\ker\pi_{*p}$ inside $T_pP$ and $A^*$ is vertical for any $A\in\lie$.

Given a connection $\omega$, the curvature $\Omega$ is described by the Mauer-Cartan equation
\be
\Omega(X,Y)=d\omega(X,Y)+\frac{1}{2}[\omega(X),\omega(Y)]
\ee
which for horizontal vectors $X,Y$ yields
\be
\omega([X,Y])=-2\Omega(X,Y)
\ee

With this notation, for any $X\in\X(M)$ there is a unique smooth {\em basic} vector field $X\n\in\X(P)$ such that $X\n_p\in Q_p$ and $\pi_*X\n=X$. In terms of these, the Lie brackets are summarized as follows.

\bp Let $X\n, Y\n$ be basic and $A^*,B^*$ vertical. Then
\be
[A^*,B^*]=[A,B]^*
\ee
\be\label{vertibra}
[X\n,B^*]=\zero
\ee
\be\label{horbra}
[X\n,Y\n]=[X,Y]\n-2\Omega(X\n,Y\n)^*
\ee

\ep

\br Suppose further that $G=\T^n=\R^n/\Gamma$, where $\Gamma$ is a lattice. This means that $\omega$ is invariant under right multiplication. Also, it follows that 
\be
\Omega=\pi^*\eta
\ee
where $\eta$ is a closed $\R^n$-valued 2-form on $M$. This can be seen by writing the connection forms downstairs at local trivializations.
\er

From this, and from the Koszul Formula, one obtains the following covariant derivative equations.
\bp\label{covdevP} Let $A,B\in\lie$ and let $X,Y\in\X(M)$. Then
\be\label{vertcov}
\nabla_{A^*}B^*=(\nabla_AB)^*,
\ee
\be\label{Stensor0}
\nabla_{X\n}Y\n=(\nabla_XY)\n+\Omega(X\n,Y\n)^*,
\ee
and 
\begin{align}
\langle \nabla_{A^*}X\n,Y\n\rangle_P&=(A,\Omega(X\n,Y\n)\rangle_G,\label{Stensor1}\\
\langle \nabla_{A^*}X\n,B^*\rangle_P&=\zero\label{Stensor2}
\end{align}

\ep
\begin{proof} Equation \eqref{horbra} and Koszul formula yield the following equalities.
\be 
\langle \nabla_{A^*}X\n,Y\n\rangle_P=\langle A,\Omega(X\n,Y\n)\rangle_G=\langle \nabla_{X\n}Y\n,A^*\rangle_P.
\ee
The first one is already the required for Equation \eqref{Stensor0} and the second one produces \eqref{Stensor1} by way of an orthonormal basis. 
The fibers being totally geodesic implies \eqref{vertcov}, and \eqref{vertibra} together with the fact that $\nabla$ is torsion-free yields \eqref{Stensor2}. 
\end{proof}
\section{Main Lemma}
\bl\label{pareqmatrix} Let $P^{n+m}$ be a principal $G^m$-bundle with connection $\omega$ over a riemannian manifold $(M^n,g_M)$ with projection map $\pi:P\rightarrow M$ and let $\Omega$ be the curvature 2-form.  Let $\xi_j$ be a local orthonormal frame on $M$, and let $A_j$ be an orthonormal frame in $G$. With respect to the orthonormal basis $\{\xi\n_i,A_j^*\}$ on $P$ along a curve $\delta$ in $P$ parallel translation is the solution to the matrix differential equation of the form 
\be
\dot Y=Q\cdot Y
\ee
with $Y(\zero)=I$, the identity matrix, and $B$ written as follows.
\be\label{pareqmatrixeq}
Q=\arraycolsep=2.5pt\def\arraystretch{2}
\left[
\begin{array}{cc}
\left(-a^k\Gamma_{kj}^i-d^k\delta_{ka}\Omega^a(\xi_j\n,\xi_{\ell}\n)\delta^{\ell i}\right)_{\overset{1\leq i\leq n}{\scriptscriptstyle1\leq j\leq n}}
& -\left(\delta_{\ell j}\Omega^{\ell}((\dot\alpha)\n, \xi_k\n)\delta^{ki}\right)_{\overset{1\leq i\leq n}{\scriptscriptstyle1\leq j\leq m}}\\ %\hline
\left(-\Omega^i((\dot\alpha)\n,\xi_j\n)\right)_{\overset{1\leq i\leq m}{\scriptscriptstyle1\leq j\leq n}}
&\left(-d^kC_{kj}^i\right)_{\overset{1\leq i\leq m}{\scriptscriptstyle1\leq j\leq m}}
\end{array}
\right]
\ee
where $\alpha$ is the projection of $\delta$,  $\dot\delta$ is given as $\dot\delta=a^i\xi_i\n+d^iA_i$, the covariant derivative in $G$ as $\nabla_{A_i}A_j=C_{ij}^kA_k$, the curvature form as $\Omega=\Omega^jA_j$, and the $\Gamma_{ij}^k$'s are the Christoffel symbols for the metric on $M$ with respect to the basis $\{\xi_i\}$: $\nabla_{\xi_i}\xi_j=\Gamma_{ij}^k\xi_k$. Caution: $\delta^{ij}$ stands for the Kronecker delta and is not related to the curve $\delta$.
\el
\begin{proof} Denote $\dot\delta$ by $T$. Let $X$ be a vector field along a curve $\delta$ in $P$, written as $X=c^i\xi_i\n+b^iA_i$
\begin{align*}
\nabla_TX&=a^ic^j\nabla_{\xi_i\n}\xi_j\n+a^i\xi_i\n(c^j)\xi_j\n+a^ib^j\nabla_{\xi_i\n}A_j^*+a^i\xi_i\n(b^j)A_j^*\\
                  &+d^ic^j\nabla_{A_i^*}\xi_j\n+d^iA_i^*(c^j)\xi_j\n+d^ib^j\nabla_{A_i^*}A_j^*+d^iA_i^*(b^j)A_j^*\\
                  &=a^ic^j(\nabla_{\xi_i}\xi_j)\n+T(c^j)\xi_j\n+a^ic^j\Omega(\xi_i\n,\xi_j\n)^*\\
                  &+a^ib^j\nabla_{\xi_i\n}A_j^*+d^ic^j\nabla_{A_i^*}\xi_j\n+T(b^j)A_j^*+d^ib^j\nabla_{A_i^*}A_j^*.
\end{align*}
once can write $a^ib^j\nabla_{\xi_i\n}A_j^*$ as follows, by \eqref{vertibra}, \eqref{Stensor1} and \eqref{Stensor2},
\begin{align*}
a^ib^j\nabla_{\xi_i\n}A_j^*&=\sum_k\langle a^ib^j\nabla_{\xi_i\n}A_j^*,\xi_k\n\rangle\xi_k\n\\
&=\sum_k a^ib^j\langle\nabla_{A_i^*}\xi_i\n,\xi_k\n\rangle\xi_k\n\\
&=\sum_k a^ib^j\langle A_j,\Omega(\xi_i\n,\xi_k\n)\rangle\xi_k\n.
\end{align*}
and  $d^ic^j\nabla_{A_i^*}\xi_j\n$ in a similar fashion,
\begin{align*}
d^ic^j\nabla_{A_i^*}\xi_j\n&=\sum_k\langle d^ic^j\nabla_{A_i^*}\xi_j\n,\xi_k\n\rangle\xi_k\n\\
%&=\sum d^ic^j\langle\nabla_{A_j^*}\xi_i\n,\xi_k\n\rangle\xi_k\n\\
&=\sum_k d^ic^j\langle A_i,\Omega(\xi_j\n,\xi_k\n)\rangle\xi_k\n,
\end{align*}
to yield
\begin{align*}
\nabla_TX&=a^ic^j(\nabla_{\xi_i}\xi_j)\n+T(c^j)\xi_j\n+c^j\Omega(a^i\xi_i\n,\xi_j\n)^*\\
                  &+\sum_k b^j\langle A_j,\Omega(a^i\xi_i\n,\xi_k\n)\rangle\xi_k\n\\
                  &+\sum_k c^j\langle d^iA_i,\Omega(\xi_j\n,\xi_k\n)\rangle\xi_k\n\\
                  &+T(b^j)A_j^*+d^ib^j\nabla_{A_i^*}A_j^*.
\end{align*}
This is further rearranged as
\begin{align*}
\nabla_TX&=\left(T(c^k)+c^ja^i\Gamma_{ij}^k+b^j\langle A_j,\Omega(a^i\xi_i\n,\xi_k\n)\rangle+c^j\langle d^iA_i,\Omega(\xi_j\n,\xi_k\n)\rangle\right)\xi_k\n\\
                 &T(b^j)A_j^*+d^ib^j\nabla_{A_i^*}A_j^*+c^j\Omega(a^i\xi_i\n,\xi_j\n)^*\\
                 &=\left(T(c^k)+c^j\left(a^i\Gamma_{ij}^k+\langle d^iA_i,\Omega(\xi_j\n,\xi_k\n)\rangle\right)+b^j\langle A_j,\Omega(a^i\xi_i\n,\xi_k\n)\rangle\right)\xi_k\n\\
                 &+\left(T(b^k)+c^j\Omega^k(a^i\xi_i\n,\xi_j\n)+b^jd^iC_{ij}^k\right)A_k^*.
\end{align*}
Therefore, $\nabla_TX$ vanishes if and only if both
$$
T(c^k)+c^j\left(a^i\Gamma_{ij}^k+\langle d^iA_i,\Omega(\xi_j\n,\xi_k\n)\rangle\right)+b^j\langle A_j,\Omega(a^i\xi_i\n,\xi_k\n)\rangle=\zero
$$
and
$$
T(b^k)+c^j\Omega^k(a^i\xi_i\n,\xi_j\n)+b^jd^iC_{ij}^k=\zero
$$
are true. Writing this as a matrix yields the claim.
\end{proof}

\section{Proof of the Homogenous Waning Theorem}\label{proof}
This section is devoted to the proof of Theorem \ref{HWT}A, stated below as Theorem \ref{homwane}. Let $(P,g_P)\stackrel{\pi}{\longrightarrow}(M,g_M)$ be principal $G$-bundle with connection and connection metric $g_P=\pi^*g_M+g_G$ for a bi-invariant metric on $G$.  Consider $\delta_n$ be a sequence of loops at $x\in P$ such that 
\be
\lim_{n\rightarrow\infty}\ell_n(\delta_n)=\zero,
\ee
where $\ell_n$ is the length with respect to the metric $g_{n}=\pi^*g_M+\varepsilon^2_ng_G$, with ${\epsilon_n}$ a sequence of positive numbers converging to zero.  Then, perhaps up to passing to a subsequence,  one would like to argue that there exists $\delta$, a vertical curve, such that
\be
\delta=\lim_{n\rightarrow\infty}\delta_n,
\ee
with respect to $g_P$. However, this is not always the case. What is true is that the sequence $\delta_n$ projects to a sequence of curves $\alpha_n$ converging to the constant curve $\alpha\equiv p=\pi(x)$. At the level of parallel translations the conclusion is stronger, as stated in the following result.

\bl\label{verthol} Let $P$ be a principal $G$-bundle with connection $\omega$ over a riemannian manifold $(M,g_M)$ with projection map $\pi:P\rightarrow M$. Let ${\epsilon_n}$ be a sequence of positive numbers converging to zero. Let $g_n$ be the connection metric on $P$ given by 
\be
g_n=\pi^*g_M+\varepsilon_n^2g_G
\ee
for a fixed bi-invariant metric $g_G$ on $G$.  Let $\delta_n$ be a sequence of curves on $P$ such that $\ell_n(\delta)$ converges to zero.  Then, up to passing to a subsequence the parallel translations $P_1^{\delta_n}$ converge,
\be
\lim_{n\rightarrow\infty}P_1^{\delta_n}=P^*,
\ee
to $P^*$, with $P$ a limit of parallel translations on $G$. 
\el

\begin{proof} Denote by $A_{n,i}=\varepsilon\inv_nA_i$. With respect to the basis $A_{n,i}$, denote by $\Omega_n=\Omega_n^jA_{n,j}$ the curvature 2-form, and let $[C_n]_{ij}^k$ be the structural constants, $\nabla_{A_{n,i}}A_{n,j}=[C_n]_{ij}^kA_{n,k}$. Let $\Gamma_{ij}^k$ be the Christoffel symbols for the metric on $M$ with respect to the basis $\{\xi_i\}$. The matrix in Equation \eqref{pareqmatrixeq} becomes
$$
Q=\arraycolsep=2.5pt\def\arraystretch{2}
\left[
\begin{array}{cc}
\left(-a^k\Gamma_{kj}^i-d^k\delta_{ka}\Omega_n^a(\xi_j\n,\xi_{\ell}\n)\delta^{\ell i}\right)_{\overset{1\leq i\leq n}{\scriptscriptstyle1\leq j\leq n}}
& -\left(\delta_{\ell j}\Omega_n^{\ell}((\dot\alpha)\n, \xi_k\n)\delta^{ki}\right)_{\overset{1\leq i\leq n}{\scriptscriptstyle1\leq j\leq m}}\\ %\hline
\left(-\Omega_n^i((\dot\alpha)\n,\xi_j\n)\right)_{\overset{1\leq i\leq m}{\scriptscriptstyle1\leq j\leq n}}
&\left(-d^k[C_n]_{kj}^i\right)_{\overset{1\leq i\leq m}{\scriptscriptstyle1\leq j\leq m}}
\end{array}
\right]
$$
with respect to any curve $\delta$ in $P$ with $\dot\delta=a^i\xi_i\n+d^iA_{n,i}$. Observe now that
$$
\Omega^j_n=\varepsilon_n\Omega^j
$$
and
$$
[C_n]_{ij}^k=\varepsilon_n\inv C_{ij}^k,
$$
wherefore the matrix becomes
$$
Q_n=\arraycolsep=2.5pt\def\arraystretch{2}
\left[
\begin{array}{cc}
\left(-a^k\Gamma_{kj}^i-\varepsilon_nd^k\delta_{ka}\Omega^a(\xi_j\n,\xi_{\ell}\n)\delta^{\ell i}\right)_{\overset{1\leq i\leq n}{\scriptscriptstyle1\leq j\leq n}}
& -\left(\varepsilon_n\delta_{\ell j}\Omega^{\ell}((\dot\alpha)\n, \xi_k\n)\delta^{ki}\right)_{\overset{1\leq i\leq n}{\scriptscriptstyle1\leq j\leq m}}\\ %\hline
\left(-\varepsilon_n\Omega^i((\dot\alpha)\n,\xi_j\n)\right)_{\overset{1\leq i\leq m}{\scriptscriptstyle1\leq j\leq n}}
&\left(-\frac{1}{\varepsilon_n}d^kC_{kj}^i\right)_{\overset{1\leq i\leq m}{\scriptscriptstyle1\leq j\leq m}}
\end{array}
\right].
$$
Finally, given a sequence of curves $\delta_n:[a,b]\rightarrow P$ with length $\ell_n(\delta)$ going to zero, if one writes $\dot\delta_n=a_n^i\xi_i\n+d_n^iA_{n,i}$ then
$$
|a^i_n|,\,|d^i_n|\leq\int_a^b\|\dot\delta(t)\|_n\,dt\rightarrow\zero
$$
and the matrices $Q_n$ satisfy the conditions of Lemma \ref{paramdepd} with $\lie$ being the standard embedding of the Lie algebra of $Hol(G)$ in the lower right of $\frak{gl}_n$.
\end{proof}

Piecing this observation together with Theorems \ref{precptBWC}, \ref{finholtype} and \ref{limsup} one gets the following theorem.
\bt[Homogeneous Waning Theorem]\label{homwane} Let $P$ be a principal $G$-bundle with connection $\omega$ over a riemannian manifold $(M,g_M)$ with projection map $\pi:P\rightarrow M$. Let ${\epsilon_n}$ be a sequence of positive numbers converging to zero. Let $g_n$ be the connection metric on $P$ given by 
\be
g_n=\pi^*g_M+\varepsilon_n^2g_G
\ee
for a fixed bi-invariant metric $g_G$ on $G$. Let $g_n^S$ be Sasaki metrics on $TP$ induced from $g_n$. Then, for any point $x\in M$ and for any point $p\in\pi\inv(x)$, the sequence of pointed spaces $(TP,g_n^S,\zero_p)$ with their projections $\pi_P$ converges in the pointed Gromov-Hausdorff sense to a space $Y$ with projection $\pi_Y:Y\rightarrow X$ with wane groups of $G_x$ given by
$$
\overline{Hol(g_G)}.
$$
\et
\begin{proof}
By Lemma \ref{verthol} one sees that the wane group $G_x$ lies inside the closure of the holonomy group of $g_G$. Now, the reverse inclusion also holds by considering sequences of vertical curves. 
\end{proof}

\bc Let $G=\T^n$, the $n$-torus group, and $P$ be a principal $G$-bundle with connection $\omega$ over a riemannian manifold $(M,g_M)$ with projection map $\pi:P\rightarrow M$. Let ${\epsilon_n}$ be a sequence of positive numbers converging to zero. Let $g_n$ be the connection metric on $P$ given by 
\be
g_n=\pi^*g_M+\varepsilon_n^2g_G
\ee
for a fixed bi-invariant metric $g_G$ on $G$. Then the wane groups are trivial,
\be
G_x=\{1\}.
\ee 
\ec

This is, for example,  the case for the Berger spheres $S^{2n-1}$ converging to $\C P^n$.

\section{Uniqueness of Parallel Translates in the Limit}\label{unipar}
This is the proof of Theorem C, stated within the section as Theorem \ref{UPT}.

\bl\label{horhol} Let $P$ be a principal $G$-bundle with connection $\omega$ over a riemannian manifold $(M,g_M)$ with projection map $\pi:P\rightarrow M$. Let ${\epsilon_n}$ be a sequence of positive numbers converging to zero. Let $g_n$ be the connection metric on $P$ given by 
\be
g_n=\pi^*g_M+\varepsilon_n^2g_G
\ee
for a fixed bi-invariant metric $g_G$ on $G$.  Let $\delta_n$ be a sequence of curves  in $P$ such that $\delta_n$ converges uniformly to a curve in $M$.  Then, up to passing to a subsequence the parallel translations $P_1^{\delta_n}$ converge,
\be
\lim_{n\rightarrow\infty}P_1^{\delta_n}=P,
\ee
to $P$, a parallel translation on $M$. 
\el

\begin{proof} Denote by $A_{n,i}=\varepsilon\inv_nA_i$. With respect to the basis $A_{n,i}$, denote by $\Omega_n=\Omega_n^jA_{n,j}$ the curvature 2-form, and let $[C_n]_{ij}^k$ be the structural constants, $\nabla_{A_{n,i}}A_{n,j}=[C_n]_{ij}^kA_{n,k}$. Let $\Gamma_{ij}^k$ be the Christoffel symbols for the metric on $M$ with respect to the basis $\{\xi_i\}$. The matrix in Equation \eqref{pareqmatrixeq} becomes
$$
Q=\arraycolsep=2.5pt\def\arraystretch{2}
\left[
\begin{array}{cc}
\left(-a^k\Gamma_{kj}^i-d^k\delta_{ka}\Omega_n^a(\xi_j\n,\xi_{\ell}\n)\delta^{\ell i}\right)_{\overset{1\leq i\leq n}{\scriptscriptstyle1\leq j\leq n}}
& -\left(\delta_{\ell j}\Omega_n^{\ell}((\dot\alpha)\n, \xi_k\n)\delta^{ki}\right)_{\overset{1\leq i\leq n}{\scriptscriptstyle1\leq j\leq m}}\\ %\hline
\left(-\Omega_n^i((\dot\alpha)\n,\xi_j\n)\right)_{\overset{1\leq i\leq m}{\scriptscriptstyle1\leq j\leq n}}
&\left(-d^k[C_n]_{kj}^i\right)_{\overset{1\leq i\leq m}{\scriptscriptstyle1\leq j\leq m}}
\end{array}
\right]
$$
with respect to any curve $\delta$ in $P$ with $\dot\delta=a^i\xi_i\n+d^iA_{n,i}$. Observe now that
$$
\Omega^j_n=\varepsilon_n\Omega^j
$$
and
$$
[C_n]_{ij}^k=\varepsilon_n\inv C_{ij}^k,
$$
wherefore the matrix becomes
$$
B_n=\arraycolsep=2.5pt\def\arraystretch{2}
\left[
\begin{array}{cc}
\left(-a^k\Gamma_{kj}^i-\varepsilon_nd^k\delta_{ka}\Omega^a(\xi_j\n,\xi_{\ell}\n)\delta^{\ell i}\right)_{\overset{1\leq i\leq n}{\scriptscriptstyle1\leq j\leq n}}
& -\left(\varepsilon_n\delta_{\ell j}\Omega^{\ell}((\dot\alpha)\n, \xi_k\n)\delta^{ki}\right)_{\overset{1\leq i\leq n}{\scriptscriptstyle1\leq j\leq m}}\\ %\hline
\left(-\varepsilon_n\Omega^i((\dot\alpha)\n,\xi_j\n)\right)_{\overset{1\leq i\leq m}{\scriptscriptstyle1\leq j\leq n}}
&\left(-\frac{1}{\varepsilon_n}d^kC_{kj}^i\right)_{\overset{1\leq i\leq m}{\scriptscriptstyle1\leq j\leq m}}
\end{array}
\right].
$$
Finally, given a sequence of curves $\delta_n:[a,b]\rightarrow P$ converging to a curve $\alpha$ in $M$, if one writes $\dot\delta_n=a_n^i\xi_i\n+d_n^iA_{n,i}$ then $|a^i_n|$ are bounded and $|d^i_n|\rightarrow\zero$ and up to passing to a subsequence the matrices $Q_n$ converge to a matrix of the form
\be\label{LPE}
Q=
\left[
\begin{array}{cc}
\left(-a^k\Gamma_{kj}^i\right)_{\overset{1\leq i\leq n}{\scriptscriptstyle1\leq j\leq n}} & \zero\\ 
\zero & \zero
\end{array}
\right],
\ee
which by Lemma \ref{paramdepd} yields the claim.
\end{proof}

\bl Let $\gamma_n$ be a sequence of horizontal curves on $TP$ over curves $\delta_n$ in $P$ such that $\gamma_n$ converges a horizontal curve in $Y$. Consider $c_n$, parallel translations along $\alpha_n=\pi\delta_n$ with $c_n(\zero)=\pi_*\gamma_n(\zero)$, $\lambda=\omega(\gamma_n(\zero))$, and let
\be
\tilde\gamma_n=c_n\n+\lambda^*
\ee
over the horizontal lift of $\alpha$ in $P$ with initial value $\pi_P\gamma_n(\zero)$.

Then, up to passing to a subsequence,
$$
\lim_{n\rightarrow\infty}\tilde\gamma_n=\lim_{n\rightarrow\infty}\gamma_n.
$$
\el
\begin{proof}
Observe that $\tilde\gamma_n$ is a solution to the limit parallel translation Equation \eqref{LPE} given in Lemma \ref{horhol} with the same initial conditions as $\gamma$. Therefore, up to passing to a subsequence, the curves $\tilde\gamma_n$ and $\gamma_n$ are as close to each other as needed to apply Arzela-Ascoli to get a convergent sequence of $\tilde\gamma_n$.
\end{proof}

%Piecing this observation together with Theorems \ref{precptBWC}, \ref{finholtype} and \ref{limsup} one gets the following theorem.
\bt \label{UPT} Let $P$ be a principal $G$-bundle with connection $\omega$ over a riemannian manifold $(M,g_M)$ with projection map $\pi:P\rightarrow M$. Let ${\epsilon_n}$ be a sequence of positive numbers converging to zero. Let $g_n$ be the connection metric on $P$ given by 
\be
g_n=\pi^*g_M+\varepsilon_n^2g_G
\ee
for a fixed bi-invariant metric $g_G$ on $G$. Let $g_n^S$ be Sasaki metrics on $TP$ induced from $g_n$. Then, for any point $x\in M$ and for any point $p\in\pi\inv(x)$, the sequence of pointed spaces $(TP,g_n^S,\zero_p)$ with their projections $\pi_P$ converges in the pointed Gromov-Hausdorff sense to a space $Y$ with projection $\pi_Y:Y\rightarrow X$ as in the Homogeneous Waning Theorem \ref{homwane}. Then parallel translation in $Y$ along curves in $M$ is unique.
\et
\begin{proof}
By Lemma \ref{horhol} and in light of Theorem \ref{limhh},  horizontal curves are limits of horizontal curves in $TP$ and these share their limits with curves of the form $c_n\n+\lambda^*$, for $c$ horizontal and $\lambda$ constant. Consider two sequences of curves of this form: $\gamma=c_n\n+\lambda^*$ and $\tilde\gamma=\tilde c_n\n+\tilde\lambda^*$. If $\lim_{n\rightarrow\infty}\tilde\gamma_n=\lim_{n\rightarrow\infty}\gamma_n$, then $\lim_{n\rightarrow\infty}\tilde c_n=\lim_{n\rightarrow\infty}c_n$ by projecting them to $TM$. Now, it must be that $\lim_{n\rightarrow\infty}\tilde\lambda_n$ and $\lim_{n\rightarrow\infty}a_n\lambda_n$ coincide for some choice of $a_n$ in the corresponding holonomy group since otherwise the initial conditions wouldn't coincide. But then $\gamma_n$ and $\tilde\gamma_n$ are uniformly getting closed to each other as the sequence progresses. This proves the claim. 
\end{proof}

\section{Local Triviality}\label{prooflt}
Theorem B in the introduction is stated here as follows.
\bt\label{LCTHM}  Let $P$ be a principal $G$-bundle with connection $\omega$ over a riemannian manifold $(M,g_M)$ with projection map $\pi:P\rightarrow M$. Let ${\varepsilon_n}$ be a sequence of positive numbers converging to zero. Let $g_n$ be the connection metric on $P$ given by 
\be
g_n=\pi^*g_M+\varepsilon_n^2g_G
\ee
for a fixed bi-invariant metric $g_G$ on $G$. Let $g_n^S$ be Sasaki metrics on $TP$ induced from $g_n$. Then, for any point $x\in M$ and for any point $p\in\pi\inv(x)$, the sequence of pointed spaces $(TP,g_n^S,\zero_p)$ with their projections $\pi_P$ converges in the pointed Gromov-Hausdorff sense to a space $Y$ with projection $\pi_Y:Y\rightarrow X$ as in the Homogeneous Waning Theorem \ref{homwane}. Then  $Y$ is locally trivial bundle over $M$ with fibers 
\be
\R^n\times(H\backslash\R^m)
\ee
\et
\begin{proof}

To produce a local trivialization, start with a small convex ball $B$ centered at $p\in M$ over which both $\pi:P\rightarrow M$ and $\pi_M:TM\rightarrow M$ are trivial.  Consider the maps $\varphi_n:TM\times\lie\rightarrow TP$ given by
 $$
 \varphi_n:(x,A)\mapsto x\n+\varepsilon_n\inv A^*
 $$
 over a fixed section $\sigma(x)\in P$.
Let $\xi_j$ be a local orthonormal frame on $M$, and let $A_j$ be an orthonormal frame in $G$.
Let $\beta$ be a geodesic in $TM$ and $\lambda$ geodesic in $\lie$. Let $\alpha=\pi_M\beta$, $\delta=\sigma\alpha$, so that $\alpha$ is the projection of $\delta$,  $\dot\delta$ is given as 
$$
\dot\delta=a^i\xi_i\n+d^iA^*_i,
$$
and denote the covariant derivative in $G$ by $\nabla_{A_i}A_j=C_{ij}^kA_k$, the curvature form by $\Omega=\Omega^jA_j$, and by $\Gamma_{ij}^k$'s  the Christoffel symbols for the metric on $M$ with respect to the basis $\{\xi_i\}$: $\nabla_{\xi_i}\xi_j=\Gamma_{ij}^k\xi_k$.
 
 Let $X$ be the vector field along a curve $\delta$ in $P$, written as $X=c^i\xi_i\n+b^iA_i^*$, where $\lambda=b^iA_i$. By the Main Lemma \ref{pareqmatrix},
\begin{align*}
\nabla_TX&=\left(T(c^k)+c^j\left(a^i\Gamma_{ij}^k+\langle d^iA_i,\Omega(\xi_j\n,\xi_k\n)\rangle\right)+b^j\langle A_j,\Omega(a^i\xi_i\n,\xi_k\n)\rangle\right)\xi_k\n\\
                 &+\left(T(b^k)+c^j\Omega^k(a^i\xi_i\n,\xi_j\n)+b^jd^iC_{ij}^k\right)A_k^*.
\end{align*}
Now, if $\gamma=X$,
\begin{align*}
\nabla_{\dot\delta}\gamma&=(\nabla_{\dot\alpha}c)\n+\sum_k\left(\big\langle\Omega(c\n ,\xi_k\n),\omega(\dot\delta)\big\rangle_G+\big\langle\Omega(\dot\alpha\n,\xi_k\n),\lambda\big\rangle_G\right)\xi_k\n\nonumber\\
 &+\left(\Omega(\dot\alpha\n,c\n)+\dot\lambda+b^jd^iC_{ij}^kA_k\right)^*\\
 &=(\nabla_{\dot\alpha}c)\n+\Omega(\dot\alpha\n,c\n)^*+\nabla_{\dot\delta}(\lambda^*)\\
                  &+\sum_k\big\langle\Omega(c\n ,\xi_k\n),\omega(\dot\delta)\big\rangle_G\xi_k\n.
\end{align*}
From which, together with the definition of $\|\dot\gamma\|$ in \eqref{sasdefeq}, it follows that $\|\dot\gamma\|$ is bounded by $\|\dot c\|+\|\dot\lambda\|+O(\varepsilon_n)$, and thus the family $\{\varphi_n\}$ is equicontinuous on $\{(X,A)|\|X\|,\|A\|\leq R\}$. The Arzel\`a-Ascoli Theorem \ref{ArAs} now gives a map $\varphi:TM\times\lie\rightarrow Y$. This map need not be injective. However, 
$$
\varphi(X,A)=\varphi(Y,B)
$$
implies that $X=Y$ and that $A\sim B\mod H$.  Indeed, $\varphi(X,A)$ can be seen as a limit of endpoints of curves of the form $$\gamma=c_n\n+\lambda^*$$ and as in theorem \ref{UPT}, two such curves must coincide up to an element in the holonomy group.  This means that one can induce a continuous function, also denoted by $\varphi$, on $TM\times(H\backslash\lie)$ where $H$ is the holonomy group of the bi-invariant metric on $G$.

Restricting to sets of the form $\{(X,A)|\|X\|,\|A\|\leq R\}$ yields that the map is an embedding. Now, since $\pi_Y\inv(B)\cap\{y|\mu(y)<r\}$ is open (cf. Remark \ref{absnorm}) in  $\pi_Y\inv(B)$ on concludes that $\varphi$ is an embedding as promised.
\end{proof}

In particular, in the case of the Berger 3-spheres, 

\bc The limit of the tangent bundles of the Berger 3-spheres, endowed with their Sasaki metrics, converge to $\R\times T\sS^2$ with the product metric.
\ec
\begin{proof}
The limit $Y$ is a rank 3 vector bundle $\xi$ over $\sS^2$. Since it contains the tangent bundle $\tau_{\sS^2}$ as a sub bundle, it must then be split (see for example the textbook by \citet{MR2045823}). 
\end{proof}

%%%%%%%%%%%%%%%%%%%%%%%%%%%%%%%%%%%%%%%%%%%%%%%%%%%%%%%%%%%%
%%%%%%%%%%%%%%%%%%%%%%%%%%%%%%%%%%%%%%%%%%%%%%%%%%%%%%%%%%%%
%%%%%%%%%%%%%%%%%%%%%%%%%%%%%%%%%%%%%%%%%%%%%%%%%%%%%%%%%%%%

\bibliographystyle{plainnat}
\bibliography{../../../ref}
\end{document}